\documentclass[12pt,a4paper]{article}
\usepackage[a4paper,margin=2cm]{geometry}

\RequirePackage{etex} 
\usepackage{amsmath}
\usepackage{amsthm}
\usepackage{subcaption}
\usepackage{amssymb}
\usepackage{enumitem}
\usepackage{graphicx}
\usepackage{nccmath}
\usepackage{setspace}
\usepackage[colorlinks=true, allcolors=blue]{hyperref}
\usepackage{enumitem}
\usepackage{multirow}
\usepackage{float}
\usepackage[dvipsnames]{xcolor}
\usepackage{tikz-network}

\newtheorem{theorem}{Theorem}
\newtheorem{proposition}[theorem]{Proposition}
\newtheorem{corollary}[theorem]{Corollary}
\newtheorem{lemma}[theorem]{Lemma}

\usepackage{thmtools}

\theoremstyle{definition}
\newtheorem{example}[theorem]{Example}

\def \mod#1{{\:({\rm mod}\ #1)}}

\def \D{\mathcal{D}}

\def \leq {\leqslant}
\def \geq {\geqslant}

\let\oldproofname=\proofname
\renewcommand{\proofname}{\textup{\textbf{\oldproofname}}}

\title{Transitive path decompositions of Cartesian products of complete graphs}

\author{Ajani De Vas Gunasekara \footnote{ajani.gunasekara@gmail.com, School of Arts and Sciences, The University of Notre Dame Australia, 140 Broadway, Sydney NSW 2007, Australia.} \qquad Alice Devillers \footnote{
alice.devillers@uwa.edu.au, Centre for the Mathematics of Symmetry and Computation,
School of Physics, Mathematics and Computing, The University of Western Australia, Crawley WA 6009, Australia.}
}

\date{}

\begin{document}

\maketitle

\begin{abstract}
An $H$-decomposition of a graph $\Gamma$  is a partition of its edge set into subgraphs isomorphic to $H$. A transitive decomposition is a special kind of $H$-decomposition that is highly symmetrical in the sense that the subgraphs (copies of $H$) are preserved and transitively permuted by a group of automorphisms of  $\Gamma$. This paper concerns transitive $H$-decompositions of the graph $K_n \Box K_n$ where $H$ is a path. When $n$ is an odd prime, we present a construction for a transitive path decomposition where the paths in the decomposition are considerably large compared to the number of vertices. Our main result supports well-known Gallai's conjecture and an extended version of Ringel's conjecture.

\end{abstract}

\let\thefootnote\relax\footnotetext{ The first author was supported by an Australian mathematical society lift-off fellowship to visit the second author. This work forms part of the ARC Discovery Grant project  DP200100080 of the second author, and was completed  when the first author was employed  at Monash University. }

\section{Introduction}

A graph $\Gamma$ is a pair $(V,E)$ where $V$ is the set of vertices and $E$ is the set of edges. We  denote the vertex set and the edge set of $\Gamma$ by $V(\Gamma)$ and $E(\Gamma)$ respectively.   All the graphs considered in this paper are undirected, unweighted simple graphs.  Let $H$ be a graph, 
an \emph{$H$-decomposition} of a graph $\Gamma = (V , E)$  is a collection $ \D$ of edge disjoint subgraphs of $\Gamma$, each isomorphic to $H$, whose edge sets partition the edge set $E$ of $\Gamma$.
The study of graph decompositions has garnered significant interest, particularly in identifying conditions under which a given graph can be decomposed into copies of a certain subgraph $H$, for example, when $H$ is a cycle, a clique or a 1-factor. 

For any graph $H$, there are three obvious necessary conditions for the existence of an $H$-decomposition of $\Gamma$. These are as follows:

\begin{lemma} [\cite{Ushio1993}]\label{L: existence of H designs}
Let $\Gamma$ be a graph. If there exists an $H$-decomposition of $\Gamma$, then
\begin{itemize}
    \item [(1)] $|V(H)| \leq |V(\Gamma)|$ or $E(\Gamma) = \emptyset$,
    \item [(2)] $|E(\Gamma)| \equiv 0 \mod{|E(H)|}$,
    \item [(3)] $\deg_\Gamma(x) \equiv 0 \mod{d}$ for each $x \in V(\Gamma)$, where $d$ is the greatest common divisor of the degrees of the vertices in $H$.
\end{itemize}
\end{lemma}

In particular, $H$-decompositions are used to solve construction problems occurring in graph theory, database systems and many related areas \cite{Ushio1993}. The general definition of $H$-decompositions was introduced by Hell and Rosa in 1972 in their work on the generalised version of the famous handcuffed prisoners problem \cite{HellRosa1972}. Notably, the solution to this problem corresponds to a fascinating construction of so-called resolvable $2$-path decomposition of $K_9$. To learn more about $H$-decompositions see \cite{AjaniThesis2022}.

An \emph{automorphism} of $\Gamma$ is a permutation of the vertex set $V$ of $\Gamma$ which leaves the edge set $E$ of $\Gamma$ invariant. Let $G$ be an automorphism group of the graph $\Gamma$. We say that the $H$-decomposition $\mathcal{D}$ of $\Gamma$ is $G$-transitive if the following two conditions hold.
\begin{enumerate}
    \item $G$ leaves $\mathcal{D}$ invariant, that is for all $H \in \mathcal{D}$ and $g \in G$, we have  $H^g\in \mathcal{D}$. 

    \item $G$ acts transitively on $\mathcal{D}$, that is 
    for any $H_1,H_2 \in \mathcal{D}$, there exists a $g \in G$ such that $H_1^g=H_2$.
   
\end{enumerate} 
If these conditions hold then we call the triple $(G, \Gamma, \mathcal{D})$ a \emph{transitive $H$-decomposition}. One of the advantages of transitivity is that one needs very little information to define or generate a decomposition.

The characterisation of graph decompositions in terms of their associated groups have attracted considerable interest over the past decades. Studying the group action of a graph is an important method for studying both the group and the graph. In other words, we study the group through its action on the graph, and we study the graph by the properties of its group. This forms a significant mathematical method and has yielded many important mathematical results. The concept of transitive decomposition of a graph has been present for a considerable period, dating back to the 1970s, and has been discussed in various forms in several publications (see \cite{LiPraeger2002, LimThesis2004, Pearce2010}). However, it is only in recent times that it has received extensive independent research attention.

In this paper, we delve into the study of transitive $H$-decompositions, focusing on cases where $H$ represents a path of a specific length. In general, path decompositions renowned for their utility in graph theory and combinatorial optimization, serve as a powerful tool for gaining valuable insights into the structural properties of graphs. They not only offer a deeper understanding of graph structure but also pave the way for the design of efficient algorithms that tackle a wide range of problems (see \cite{Anshelevich Zhang2004, KrahnBinderKonig2014}).

Path decompositions of various types of graphs have been studied extensively, including complete graphs \cite{Tarsi1983}, bipartite graphs \cite{ChuFanZhou2021}, and digraphs (such as tournaments \cite{Girao2023}). The problem of decomposing digraphs into paths was initially investigated by Alspach and Pullman \cite{AlspachPullman1974}, who established bounds on the minimum number of paths necessary for such decompositions.  However, there is limited literature on path decompositions of graph products. In \cite{JeevadossMuthusamy2016}, Jeevadoss and Muthusamy discuss path and cycle decompositions of product graphs, particularly focusing on cases where the path and cycle lengths are small.

A group \(G\) acting on the edges \(E(\Gamma)\) of a graph \(\Gamma\) is said to be \emph{semiregular} if no element of \(G\) other than the identity fixes an edge of \(\Gamma\). This means that each non-identity element of \(G\) maps every edge to a different edge.
We give a construction for a transitive path decomposition of $\Gamma = K_n \Box K_n$, when $n$ is an odd prime and the group $G$ acting on $\Gamma$ is semiregular on $E(\Gamma)$. Note that the paths in the decomposition are considerably long.

\begin{theorem} \label{T: main theorem}
A transitive $n(n-1)$-path decomposition of $K_n \Box K_n$ exists for all odd primes $n$.
\end{theorem}
Note this implies in particular that for any arbitrarily large number $t$, there is a graph with a transitive path decomposition with path length at least $t$.

Gallai's well-known conjecture (see \cite{Bonamy2019}) states that, every connected graph on $m$ vertices can be decomposed into at most $\frac{m+1}{2}$ paths. In \cite{Lovasz1968}, Lovász proved that every graph on $m$ vertices can be decomposed into at most $\frac{m}{2}$ paths and cycles. Moreover they showed that every graph on $m$ vertices with odd vertex degrees can be decomposed into $\frac{m}{2}$ paths.
There are many partial results available towards Gallai's conjecture, including asymptotic results (see \cite{Pyber1996}).
In our construction proving Theorem~\ref{T: main theorem},  $K_n \Box K_n$ has $m = n^2$ vertices (where $n$ is an odd prime), and the decomposition has $n=\sqrt{m}$ paths, since $|E(K_n \Box K_n)| = n^2(n-1)$. This validates Gallai's conjecture in this case.

A generalization of Ringel's conjecture, attributed to Häggkvist (see \cite{Hiiggkvist1989}), proposes that every $2b$-regular graph can be decomposed into trees with $b$ edges. Theorem~\ref{T: main theorem} validates this for $K_n \Box K_n$ when $n$ is an odd prime. Indeed $K_n \Box K_n$ is a $2(n-1)$-regular graph and the constructed paths can be decomposed into $n$ distinct $(n-1)$-paths. Therefore, our main result states that when $n$ is an odd prime $K_n \Box K_n$ admits a decomposition into $n^2$  paths with $n-1$ edges.

In \cite[Chapters 6--8]{PearceThesis2007}, Pearce has studied rank $3$ transitive decompositions of $K_n \Box K_n$  in general. In particular, he characterises possible subgraphs $H$ in $K_2 \Box K_2$ which admits a $G$-transitive $H$-decomposition when $G$ is a certain group.

Tarsi proved that the obvious condition $|E(\lambda K_n)| \equiv 0 \mod {\ell}, n \geq \ell +1 $ which is necessary for the existence of a decomposition of a complete multigraph $\lambda K_n$ into paths of length $\ell$ is also sufficient \cite{Tarsi1983}. Among the other things, in this paper, we show that the necessary condition $|E(K_n \Box K_n)| \equiv 0 \mod {n(n-1)}$ is also sufficient for the existence of a  decomposition into paths of length $n(n-1)$ of $K_n \Box K_n$ for each odd prime $n$.

\section{Preliminaries}

This section will establish the foundation for our work. 
The following theorem, which is of independent interest in its own right, establishes a crucial necessary condition for the existence of a transitive $H$-decomposition of a graph $\Gamma$.

\begin{theorem} \label{T: transitive decomposition}
Let $\Gamma$ be a graph and $H$ be a subgraph of $\Gamma$. Suppose that $G\leq \mathrm{Aut}(\Gamma)$ is semiregular on the edges of $\Gamma$ and $H$ contains exactly one edge from each edge orbit of $\Gamma$ under $G$. Then $H^G$ is a $G$-transitive $H$-decomposition of $\Gamma$ of size $|G|$.
\end{theorem}

\begin{proof}
Assume that $G\leq \mathrm{Aut}(\Gamma)$ is semiregular on the edges of $\Gamma$  and $H$ contains exactly one edge from each edge orbit of $\Gamma$ under $G$. Let $\mathcal{D} =H^G=\{H^g|g \in G\}$. We will first prove that $\mathcal{D}$ is a  decomposition of $\Gamma$.

Suppose that there exists an edge $e \in E(\Gamma)$ such that $e \in E(H^{g_1}) \cap E(H^{g_2})$ for some distinct $g_1, g_2 \in G$. 
Then $e^{g_1^{-1}}$ and $e^{g_2^{-1}}$ are edges in $E(H)$ that lie in the same edge orbit under $G$. Since $H$ contains exactly one edge from each edge orbit under $G$, we have $e^{g_1^{-1}}=e^{g_2^{-1}}$. This implies $e^{g_1^{-1}g_2} = e$. Since $G$ is semiregular on $E(\Gamma)$, it follows that $g_1^{-1}g_2 = I$ (the identity element of $G$) and hence $g_1 = g_2$, which is a contradiction. Thus, each edge of $\Gamma$ belongs to at most one element of $H^G$. 

Let $e \in E(\Gamma)$. Since $H$ contains exactly one edge from each edge orbit, there exists some $e_1 \in E(H)$ such that $e \in e_1^G$. This implies $e \in H^g$ for some $g \in G$. Therefore each edge of $\Gamma$ belongs to at least one element of $H^G$.

Thus each edge of $\Gamma$ belongs to exactly one element of $H^G$, and each subgraph is isomorphic to $H$ by construction. In other words $\mathcal{D}$ is an $H$-decomposition of $\Gamma$. 
Due to the nature of $G$ and $\mathcal{D}$, $G$ leaves $\mathcal{D}$ invariant and $G$ acts transitively on $\mathcal{D}$, so $\mathcal{D}$ is a $G$-transitive $H$-decomposition of $\Gamma$.

Since $H$ contains exactly one edge from each edge orbit, the stabiliser $G_H$ of $H$ fixes all edges of $H$ and thus is trivial, as $G$ is semiregular on edges. Hence $|G|=|H^G|\cdot |G_H|=|H^G|=|\mathcal{D}|$.
  
\end{proof}

The above theorem is applied in the following example to illustrate the process of finding a $G$-transitive $H$-decomposition of $K_9$, where $H$ comprises four edge-disjoint triangles (copies of $K_3)$, and $G\leq \mathrm{Aut}(K_9)$. Moreover, this particular decomposition also yields a $K_3$-decomposition of $K_9$.

\begin{example}
Let $\Gamma = K_9$ and we denote the vertex set of $\Gamma$ by $\{1,2,3, \ldots, 9\}$. Let $G = \langle (1,4,7)(2,5,8)(3,6,9) \rangle$ of order $3$. We have listed the edge orbits of $K_9$ under $G$ and highlighted the edges of $H$.

\begin{tabular}[l]{lll}
Orbit 1: & $[\textbf{\{1,4\}}, \{4,7\}, \{7,1\}]$ & Orbit 2:  $[\{2,5\}, \{5,8\}, \textbf{\{8,2\}}]$ \\
Orbit 3: & $[\{3,6\}, \{6,9\}, \textbf{\{9,3\}}]$ & Orbit 4:  $[\{1,2\}, \textbf{\{4,5\}}, \{7,8\}]$ \\
Orbit 5: & $[\textbf{\{1,5\}}, \{4,8\}, \{7,2\}]$ & Orbit 6:  $[\{1,6\}, \{4,9\}, \textbf{\{7,3\}}]$ \\
Orbit 7: & $[\{2,3\}, \textbf{\{5,6\}}, \{8,9\}]$ & Orbit 8:  $[\{2,4\}, \textbf{\{5,7\}}, \{8,1\}]$ \\
Orbit 9:  &$[\textbf{\{2,6\}}, \{5,9\}, \{8,3\}]$ & Orbit 10:  $[\textbf{\{3,7\}}, \{6,1\}, \{9,4\}]$ \\
Orbit 11: & $[\{3,5\}, \textbf{\{6,8\}}, \{9,2\}]$ & Orbit 12:  $[\textbf{\{3,9\}}, \{6,3\}, \{9,6\}]$ \\
\end{tabular}

\begin{figure}[ht]
  \begin{minipage}[t]{0.5\linewidth}
    \centering
    \begin{tikzpicture} [scale=1]
\foreach \i in {1,...,9}
    \node[circle,fill=black,inner sep=1.5pt, ] (\i) at ({360/9 * (\i - 1)}:2.5) {};

  \foreach \i/\label in {1/1, 2/2, 3/3, 4/4, 5/5}
    \node[circle,fill=black,inner sep=1.5pt,label=above:\label] at (\i) {};
  
  \foreach \i/\label in {6/6, 7/7, 8/8, 9/9}
    \node[circle,fill=black,inner sep=1.5pt,label=below:\label] at (\i) {};

  \draw[color = black, line width = 1.5pt] (1) -- (4);
  \draw[color = black, line width = 1.5pt] (4) -- (5);
  \draw[color = black, line width = 1.5pt] (1) -- (5);

  \draw[color = black, line width = 1.5pt] (5) -- (6);
  \draw[color = black, line width = 1.5pt] (6) -- (7);
  \draw[color = black, line width = 1.5pt] (7) -- (5);

  \draw[color = black, line width = 1.5pt] (2) -- (6);
  \draw[color = black, line width = 1.5pt] (2) -- (8);
  \draw[color = black, line width = 1.5pt] (6) -- (8);

  \draw[color = black, line width = 1.5pt] (7) -- (9);
  \draw[color = black, line width = 1.5pt] (9) -- (3);
  \draw[color = black, line width = 1.5pt] (7) -- (3);

    \end{tikzpicture}
    \caption{$H$}
    \label{fig:first}
  \end{minipage}%
  \begin{minipage}[t]{0.5\linewidth}
    \centering
    \begin{tikzpicture}[scale=1]
\foreach \i in {1,...,9}
    \node[circle,fill=black,inner sep=1.5pt, ] (\i) at ({360/9 * (\i - 1)}:2.5) {};

  \foreach \i/\label in {1/1, 2/2, 3/3, 4/4, 5/5}
    \node[circle,fill=black,inner sep=1.5pt,label=above:\label] at (\i) {};
  
  \foreach \i/\label in {6/6, 7/7, 8/8, 9/9}
    \node[circle,fill=black,inner sep=1.5pt,label=below:\label] at (\i) {};

  \draw[color = black, line width = 1.5pt] (1) -- (4);
  \draw[color = black, line width = 1.5pt] (4) -- (5);
  \draw[color = black, line width = 1.5pt] (1) -- (5);

  \draw[color = black, line width = 1.5pt] (5) -- (6);
  \draw[color = black, line width = 1.5pt] (6) -- (7);
  \draw[color = black, line width = 1.5pt] (7) -- (5);

  \draw[color = black, line width = 1.5pt] (2) -- (6);
  \draw[color = black, line width = 1.5pt] (2) -- (8);
  \draw[color = black, line width = 1.5pt] (6) -- (8);

  \draw[color = black, line width = 1.5pt] (7) -- (9);
  \draw[color = black, line width = 1.5pt] (9) -- (3);
  \draw[color = black, line width = 1.5pt] (7) -- (3);

  \draw[color=red, line width=1.5pt] (4) -- (7);
  \draw[color=red, line width=1.5pt] (7) -- (8);
  \draw[color=red, line width=1.5pt] (4) -- (8);

  \draw[color=red, line width=1.5pt] (8) -- (9);
  \draw[color=red, line width=1.5pt] (9) -- (1);
  \draw[color=red, line width=1.5pt] (1) -- (8);

  \draw[color=red, line width=1.5pt] (5) -- (9);
  \draw[color=red, line width=1.5pt] (5) -- (2);
  \draw[color=red, line width=1.5pt] (9) -- (2);

  \draw[color=red, line width=1.5pt] (1) -- (3);
  \draw[color=red, line width=1.5pt] (3) -- (6);
  \draw[color=red, line width=1.5pt] (1) -- (6);

   \draw[color=blue, line width=1.5pt] (4) -- (3);
  \draw[color=blue, line width=1.5pt] (2) -- (3);
  \draw[color=blue, line width=1.5pt] (4) -- (2);

  \draw[color=blue, line width=1.5pt] (2) -- (7);
  \draw[color=blue, line width=1.5pt] (2) -- (1);
  \draw[color=blue, line width=1.5pt] (1) -- (7);

  \draw[color=blue, line width=1.5pt] (6) -- (9);
  \draw[color=blue, line width=1.5pt] (4) -- (6);
  \draw[color=blue, line width=1.5pt] (9) -- (4);

  \draw[color=blue, line width=1.5pt] (5) -- (3);
  \draw[color=blue, line width=1.5pt] (3) -- (8);
  \draw[color=blue, line width=1.5pt] (5) -- (8);

    \end{tikzpicture}
    \caption{$G$-transitive $H$-decomposition of $K_9$}
    \label{fig:second}
  \end{minipage}
\end{figure}

Figure~\ref{fig:first} is the subgraph $H$ of $K_9$ consisting of a union of edge disjoint triangles $\{1,4,5\}$, $\{2,6,8\}$,$ \{3,7,9\}$ and $\{5,6,7\}$. Then $H^G$ gives a $G$-transitive $H$-decomposition as in Figure~\ref{fig:second}. The decomposition of $H$ into four triangles thus induces a $K_3$-decomposition of $K_9$ (but not $G$-transitive). However, this $K_3$-decomposition is transitive under different groups. One such  example  is when $G = \langle  (1, 5, 6, 8)(2, 9, 4, 7), (1, 9, 6, 7)(2, 8, 4, 5), (1, 6)(2, 4)(5, 8)(7, 9), \\ (1, 5, 4)(2, 8, 6)(3, 9, 7), (1, 9, 8)(2, 4, 3)(5, 7, 6)  \rangle$.
\end{example}

A \emph{walk} $W$ in graph $\Gamma$ is a sequence $v_0, e_1,v_1, \ldots, v_{k-1}, e_k, v_k $ of $v_i \in V(\Gamma)$ and $e_i \in E(\Gamma)$ such that for $1 \leq i \leq k$, $e_i = \{v_{i-1},v_i\}$ and denoted as $W = v_0v_1 \ldots v_k$. We define $V(W) = \{v_0,v_1, \ldots, v_k\}$ and $E(W) = \{ e_1, e_2, \ldots, e_k\}$. Since our focus is on the subgraph of $\Gamma$ formed by the edges $E(W)$ and the vertices $V(W)$, the order of the sequence is irrelevant.
Note that, in a walk edges can be repeated. The length of a walk is its number of edges (including repeated edges). A walk is known as an $\ell$-walk if its length is $\ell$.
If all the $v_i$s are distinct and hence all the $e_i$s are distinct, then the walk is known as a \emph{path} $P$ in $\Gamma$. We denote a path of length $\ell$ by $P_{\ell+1}$.

Let $\Gamma = K_n \Box K_m$ be the \emph{Cartesian product} of the complete graphs $K_n$ and $K_m$. The graph $\Gamma$ may be viewed as a $2$-dimensional `grid' consisting of `horizontal' copies of $K_m$ and `vertical' copies of $K_n$. Let $V(\Gamma)$ and $E(\Gamma)$ be the vertex set and edge set of $\Gamma$ respectively.  Let $\mathbb{Z}_n$ be the additive group of integers modulo $n$. 
We label the vertices of $\Gamma$ as ordered pairs $(a,b) \in \mathbb{Z}_n \times \mathbb{Z}_m$. Then we define the edges of $\Gamma$ as follows: for any $(a,b), (c,d) \in V(\Gamma)$, $\{(a,b),(c,d)\} \in E(\Gamma)$ whenever $a = c$ or $b = d$. 
When $a=c$, that is $(a,b)-(c,d) \in \{0\} \times \mathbb{Z}_m^* $ we call the edge $\{(a,b),(c,d)\}$ a horizontal edge and when $b = d$, that is $(a,b)-(c,d) \in \mathbb{Z}_n^* \times \{0\}$ we call the edge $\{(a,b),(c,d)\}$ a vertical edge. Unless explicitly mentioned, $\Gamma$ will denote the Cartesian product of two complete graphs throughout the remainder of this paper.

\begin{example}
$\Gamma = K_3 \Box K_4$ where each horizontal line induces a $K_4$ and each vertical line induces a $K_3$.

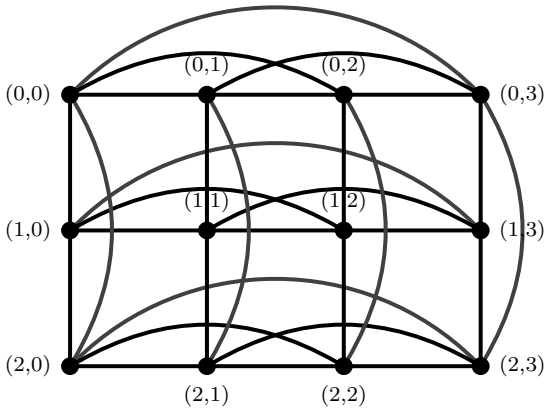
\begin{tikzpicture}[scale=0.9]

\Vertex[x=1,y=6,label = $\text{(0,0)}$, size = 0.2, color = black, position = left]{A}
\Vertex[x=3,y=6,label = $\text{(0,1)}$, size = 0.2, color = black, position = above]{B}
\Vertex[x=5,y=6,label = $\text{(0,2)}$, size = 0.2, color = black , position = above]{C}
\Vertex[x=7,y=6,label = $\text{(0,3)}$, size = 0.2, color = black , position = right]{D}

\Vertex[x=1,y = 4, label = $\text{(1,0)}$, size = 0.2, color = black, position = left]{E}
\Vertex[x=3, y= 4, label = $\text{(1,1)}$, size = 0.2, color = black, position = above]{F}
\Vertex[x=5,y=4, label = $\text{(1,2)}$, size = 0.2, color = black, position = above]{G}
\Vertex[x=7,y=4, label = $\text{(1,3)}$, size = 0.2, color = black, position = right]{H}

\Vertex[x=1,y=2, label = $\text{(2,0)}$, size = 0.2, color = black, position = left]{I}
\Vertex[x=3,y=2, label = $\text{(2,1)}$, size = 0.2, color = black, position = below]{J}
\Vertex[x=5,y=2, label = $\text{(2,2)}$, size = 0.2, color = black, position = below]{K}
\Vertex[x=7,y=2, label = $\text{(2,3)}$, size = 0.2, color = black, position = right]{L}

\Edge[color = black](A)(B)
\Edge[bend = 45](A)(D)
\Edge[color = black, bend = 30](A)(C)
\Edge[color = black](B)(C)
\Edge[color = black, bend = 30](B)(D)
\Edge[color = black](C)(D)

\Edge[color = black](E)(F)
\Edge[bend = 45](E)(H)
\Edge[color = black, bend = 30](E)(G)
\Edge[color = black](F)(G)
\Edge[color = black, bend = 30](F)(H)
\Edge[color = black](G)(H)

\Edge[color = black](I)(J)
\Edge[bend = 45](I)(L)
\Edge[color = black, bend = 30](I)(K)
\Edge[color = black](J)(K)
\Edge[color = black, bend = 30](J)(L)
\Edge[color = black](K)(L)

\Edge[color = black](A)(E)
\Edge[bend = 30](A)(I)
\Edge[color = black](E)(I)

\Edge[color = black](B)(F)
\Edge[bend = 30](B)(J)
\Edge[color = black](F)(J)

\Edge[color = black](C)(G)
\Edge[bend = 30](C)(K)
\Edge[color = black](G)(K)

\Edge[color = black](D)(H)
\Edge[bend = 30](D)(L)
\Edge[color = black](H)(L)

\end{tikzpicture}    
\end{example}

Let $W = v_0v_1v_2 \ldots v_{\ell-1}v_\ell$ be an $\ell$-walk  in $\Gamma$ and $V(W)$ and $E(W)$ be the vertex set and edge set of $W$ respectively.  
For a given walk $W$ with fixed $v_0$, we can define an array $\overrightarrow{a} = [a_1,a_2, \ldots, a_{\ell}]$ such that $a_i = v_i - v_{i-1}$ for $1 \leq i \leq \ell$ and $|\overrightarrow{a}| = \ell$. This implies $a_i \in \{(c,0),(0,c')|c\in \mathbb{Z}_n^*  \text{ and } c'\in \mathbb{Z}_m^*\}$. Conversely, a given $v_0$ and array $\overrightarrow{a}$ determines $W$ since $v_i = v_0 + \sum_{g= 1}^i a_g$. 
 We write $W = W(v_0, \overrightarrow{a})$. 
Note that we can choose $v_\ell$ as the starting vertex of $W$, in this case the array would be the negative of the reverse of $\overrightarrow{a}$, and we consider that the same walk. 
If $W(v_0, \overrightarrow{a})$ is a path, then we denote it by $P(v_0, \overrightarrow{a})$. The subsequent lemma provides the necessary and sufficient conditions
that determine whether a walk associated with an array $\overrightarrow{a}$ in $\Gamma$ is a path.

\vspace{0.5 cm}
\begin{lemma} \label{L: walk is a path}
An $\ell$-walk $W = W(v_0, \overrightarrow{a})$ is an $\ell$-path if and only if $\sum_{g= i+1}^j a_g \neq (0,0)$ for all $0 \leq i < j \leq \ell$.

\end{lemma}

\begin{proof}
An $\ell$-walk is an $\ell$-path if and only if $v_i \neq v_j$ whenever $i <j$ for all $v_i,v_j \in V(W)$. 
That is an $\ell$-walk is an $\ell$-path if and only if  $ v_0 + \sum_{g= 1}^i a_g \neq v_0 + \sum_{g= 1}^j a_g$, simplified to   $\sum_{g= i+1}^j a_g \neq (0,0)$, for all $0 \leq i < j \leq \ell$.

\end{proof}

Throughout this paper, the group $G$ is always defined as follows unless stated otherwise.  Consider the cyclic permutation of $V(\Gamma)$ defined by  $c(a,b) = (a+1,b)$ for all $(a,b) \in \mathbb{Z}_n \times \mathbb{Z}_m$. Let $G = \langle c \rangle$ which permutes the rows cyclically. This implies that the order of $G$ is $n$.

Note that all horizontal edges are mapped to horizontal edges and all vertical edges are mapped to vertical edges under $G$. Each horizontal edge orbit under $G$ can be denoted as the set $\{\{(a,b_1), (a,b_2)\} \in E(\Gamma)| a \in \mathbb{Z}_n\}$ for fixed $b_1 \neq b_2\in \mathbb{Z}_m$, while each vertical edge orbit under $G$ can be denoted as the set $\{\{(a,b), (a+t,b) \} \in E(\Gamma)| a \in \mathbb{Z}_n\}$ for fixed $ t\in\mathbb{Z}_n^*$ and $b\in \mathbb{Z}_m$. Note that the vertical orbits for $t$ and $n-t$ are the same orbit so we may assume that $0<t \leq \lfloor\frac{n-1}2\rfloor$. 
The goal of the rest of this section will be to characterise the action of the group $G$ on the graph $K_n \Box K_m$.

\begin{lemma} \label{L: number of edge orbits}
Let $n$ be an odd integer and $\Gamma = K_n \Box K_m$.
\begin{itemize}
    \item [(a)] There are $\binom{m}{2}$ horizontal edge orbits of size $n$ under $G$.
    \item [(b)] There are $\frac{m(n-1)}{2}$ vertical edge orbits of size $n$ under $G$.
\end{itemize}
\end{lemma}

\begin{proof}
  (a)   Due to the notion of horizontal edge orbits, there are $n$ edges in each horizontal edge orbit. Moreover, recall that each horizontal line of $\Gamma$ is a copy of $K_m$. Therefore, there are $\binom{m}{2}$ horizontal edge orbits of size $n$ under $G$.

(b)   We know that $G$ permutes the rows cyclically and all vertical edges are mapped to vertical edges under $G$. Consider a vertical edge $\{(a,b), (a+t,b)\}$ for some fixed $t\in\mathbb{Z}_n^*$ and $b\in \mathbb{Z}_m$. The vertical edge orbit under $G$ is given by $\{\{(a,b), (a+t,b) \} \in E(\Gamma)| a \in \mathbb{Z}_n\}$. Since the vertical orbits for $t$ and $n-t$ are the same orbit, there are $\frac{n-1}{2}$ possible values of $t$ noting that $n$ is odd, and moreover, each vertical edge orbit contains $n$ edges since $a \in \mathbb{Z}_n$. Therefore, there are $\frac{m(n-1)}{2}$ vertical edge orbits, each of size $n$ since $\Gamma$ has $m$ vertical lines (columns).
\end{proof}

By the orbit-stabilizer theorem, it follows that the stabilizer $G_e$ of $e$ under the group $G$ is trivial thus we have the following corollary.
\begin{corollary}\label{C: semiregular}
 Let $n$ be an odd integer and $\Gamma = K_n \Box K_m$. Then $G$ acts semiregularly on $E(\Gamma)$.
\end{corollary}

We denote the first coordinate and second coordinate of an ordered pair $v$ by $v_1$ and $v_2$ respectively. If $v = (a,b)$, then $v_1 =  a$ and $v_2 =  b$.

\begin{lemma} \label{L: edge orbit}
Any two distinct edges $\{v,w\}$ and $\{v', w'\}$ of $E(\Gamma)$ are in the same edge orbit under $G$ if and only if
\begin{enumerate}
    \item [(a)] $v_2 = v'_2$ and $w-v = w'-v'$ or,
    \item [(b)] $v_2 = w'_2$ and $w-v = v'-w'$.
\end{enumerate}

\end{lemma}

\begin{proof}
First, suppose that $\{v,w\}$ and $\{v', w'\}$ are horizontal edges. Let  $v = (a,v_2), w = (a, w_2), v' = (a',v'_2), w' = (a',w'_2)$ for some $a,a' \in \mathbb{Z}_n$ and $v_2,w_2,v'_2,w'_2 \in \mathbb{Z}_m$. Thus, according to the characterisation of horizontal edge orbits under $G$,  $\{v,w\}$ and $\{v', w'\}$ are in the same horizontal edge orbit if and only if $v_2 = v'_2$ and $w_2 = w'_2$, or $v_2 = w'_2$ and $w_2=v'_2$. This is equivalent to (a) or (b), respectively.

Now suppose that $\{v,w\}$ and $\{v', w'\}$ are vertical edges. Let  $v = (v_1,v_2), w = (v_1+t, v_2), v' = (v'_1,v'_2), w' = (v'_1+r,v'_2)$ for some $t,r \in \mathbb{Z}_n^*, v_1,v_1' \in \mathbb{Z}_n$ and $v_2,v'_2 \in \mathbb{Z}_m$ (note $v_2=w_2$ and $v'_2=w'_2$ in this case). Thus, according to the characterisation of vertical edge orbits under $G$, $\{v,w\}$ and $\{v', w'\}$ are in the same vertical edge orbit if and only if $v_2 = v'_2$ and $r = t$, or $v_2 = v'_2$ and $r = n-t=-t$ because the vertical orbits for $t$ and $n-t$ are the same orbit. This is equivalent to (a) or (b), respectively.
\end{proof}

The next lemma presents both the necessary and sufficient conditions that guarantee a path does not contain more than one edge from each edge orbit under the action of the group $G$.

\begin{lemma}\label{L: edge orbit general}
 Consider an $\ell$-path $P = P(v_0, \overrightarrow{a})$ in $\Gamma$. For $0 \leq i <j < \ell$ we define the following conditions:
\begin{enumerate} [label=(\alph*)]
    \item  $(\sum_{g=i+1}^j a_g)_2 \neq 0$ or $a_{i+1} \neq a_{j+1}$,
    \item $(\sum_{g=i+1}^{j+1} a_g)_2 \neq 0$ or $a_{i+1} \neq -a_{j+1}$.
\end{enumerate}
Then $P$ contains at most one edge from each edge orbit under $G$ if and only if (a) and (b) hold for all $i,j \in \{0,1, \ldots, \ell-1\}$ with $i <j$.

\end{lemma}

\begin{proof}
We will prove that $P$ contains more than one edge from an edge orbit under $G$ if and only if $\lnot$(a) $a_{i+1} = a_{j+1}$ and $(\sum_{g=i+1}^j a_g)_2 = 0$   or
$\lnot$(b) $a_{i+1} =-a_{j+1}$ \text{ and } $ (\sum_{g=i+1}^{j+1} a_g)_2 = 0$ hold for some $i,j$ such that $0 \leq i < j \leq \ell -1$.
Let $\{v_i,v_{i+1} \}\text{ and } \{v_j,v_{j+1}\} \in E(P)$. Note that, $v_{i+1} = v_i+a_{i+1} = v_0 + \sum_{g=1}^{i+1}a_g$ and $v_{j+1} = v_j+a_{j+1} = v_0 + \sum_{g=1}^{j+1}a_g$ by the definition of $\overrightarrow{a}$.
Then by Lemma~\ref{L: edge orbit}, $\{v_i,v_{i+1} \}\text{ and } \{v_j,v_{j+1}\}$ are in the same edge orbit under $G$, if and only if $(v_i)_2 = (v_j)_2$ and $a_{i+1} = a_{j+1}$ or $(v_i)_2 = (v_{j+1})_2$ and $a_{i+1} = -a_{j+1}$.

Note that, $(v_i)_{2} = (v_j)_{2}$ is equivalent to $(\sum_{g=i+1}^j a_g)_2 = 0$ since $(v_i)_2=(v_0 + \sum_{g=1}^{i}a_g)_2$, $(v_j)_2 = (v_0 + \sum_{g=1}^{j}a_g)_2$ and $i<j$. Similarly, $(v_i)_2 = (v_{j+1})_2$ is equivalent to  $(\sum_{g=i+1}^{j+1} a_g)_2 = 0$.

\end{proof}

\section{Staircase array and main results}
We will now focus on $\Gamma = K_n \Box K_n$. The goal of this section is to provide a proof of our main theorem. We aim to construct a path decomposition of $K_n \Box K_n$ that is $G$-transitive. The following example illustrates a $6$-path in $K_3 \Box K_3$. Recall that $G = \langle c \rangle$ where $c(a,b) = (a+1,b)$ for all $(a,b) \in \mathbb{Z}_3 \times \mathbb{Z}_3$.

\begin{example}
A $6$-path $P=P((0,0), \overrightarrow{a})$ in $K_3 \Box K_3$   with $\overrightarrow{a} = [(0,1),(1,0),(0,1),(1,0),(0,1),(2,0)]$. 
\begin{figure}[ht]
    
\begin{tikzpicture} [scale= 0.75]

\Vertex[x=1,y=8,label = $\text{(0,0)}$, size = 0.2, color = black, position = left]{A}
\Vertex[x=4,y=8,label = $\text{(0,1)}$, size = 0.2, color = black, position = above]{B}
\Vertex[x=7,y=8,label = $\text{(0,2)}$, size = 0.2, color = black , position = right]{C}
\Vertex[x=1,y = 5, label = $\text{(1,0)}$, size = 0.2, color = black, position = left]{D}
\Vertex[x=4, y= 5, label = $\text{(1,1)}$, size = 0.2, color = black, position = above]{E}
\Vertex[x=7,y=5, label = $\text{(1,2)}$, size = 0.2, color = black, position = right]{F}
\Vertex[x=1,y=2, label = $\text{(2,0)}$, size = 0.2, color = black, position = left]{G}
\Vertex[x=4,y=2, label = $\text{(2,1)}$, size = 0.2, color = black, position = below]{H}
\Vertex[x=7,y=2, label = $\text{(2,2)}$, size = 0.2, color = black, position = right]{I}

\draw[color = red , line width = 1.5 pt] (A) -- (B);
\Edge[bend = 45](A)(C)
\Edge[color = black](B)(C)

\Edge[color = black](D)(E)
\Edge[bend = 45, color = black](D)(F)
\draw[color = red , line width = 1.5 pt] (E) -- (F);

\Edge[color = black](G)(H)
\Edge[bend=45, color= red](G)(I)
\Edge[color = black](H)(I)

\Edge[color = black](A)(D)
\Edge[bend = 45, color = black](A)(G)
\draw[color = red , line width = 1.5 pt](D) -- (G);

\draw[color = red , line width = 1.5 pt] (B) -- (E);
\Edge[bend = 45, color = black](B)(H)
\Edge[color = black](E)(H)

\Edge[color = black](C)(F)
\Edge[bend = 45, color = black](C)(I)
\draw[color = red , line width = 1.5 pt] (F) -- (I);

\end{tikzpicture}

\captionsetup{justification=raggedright,singlelinecheck=false}
\caption{The $6$-path is represented by the red edges.}
\label{fig: third}
\end{figure}
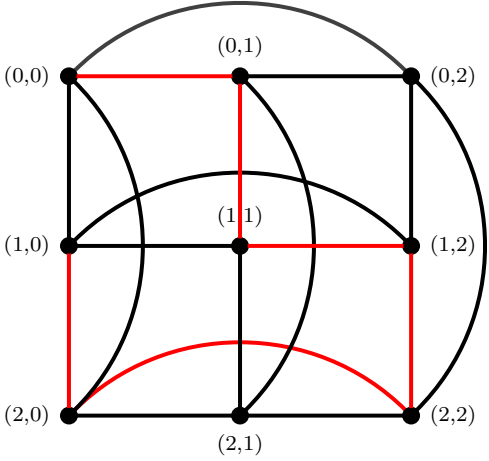

Note that $K_3 \Box K_3$ admits $3$ horizontal edge orbits and $3$ vertical edge orbits under $G$ and that $P$ contains exactly one edge from each edge orbit. By Corollary~\ref{C: semiregular} we have that $G$ is semiregular on edges. Then by Theorem~\ref{T: transitive decomposition}, $P^G$ is a $G$-transitive decomposition of $K_3 \Box K_3$  consisting of three paths, namely $P=P((0,0), \overrightarrow{a})$, $P'=P^c=P((1,0), \overrightarrow{a})$, $P''=P^{c^2}=P((2,0), \overrightarrow{a})$.

\end{example}
We will now generalize this example. 
\vspace{0.5 cm}

Recall that, $|\overrightarrow{a}|$ denotes the length of the array. For an $\ell$-walk $W = W(v_0, \overrightarrow{a})$, we have $\ell = |E(W)| = |\overrightarrow{a}| $ and $a_{i} = v_{i}-v_{i-1}$ for all $\{v_{i-1}, v_{i}\} \in E(W)$ where $0 < i \leq \ell$. 
For an odd integer $n$, we refer to an array $\overrightarrow{a}$ in the form of

\begin{multline*}
[(0,1),(1,0), \ldots, (0,1),(2,0),(0,3),(3,0),\ldots,(0,3),(4,0), \\ \ldots, (0,n-2),(n-2,0)\ldots,(0,n-2),(n-1,0)]
\end{multline*}
as the \emph{staircase array}. This array is composed of a concatenation of subarrays of length $2n$, called \emph{stretches} of the form $[(0,c),(c,0),(0,c)\dots, (c,0),(0,c),(c+1,0)]$ where $c$ takes the values $2k - 1$ for $1 \leq k \leq (n-1)/2$ noting that $n$ is odd. Therefore, a staircase array $\overrightarrow{a}$ has $\frac{n-1}{2}$ such stretches. Thus $|\overrightarrow{a}| = 2n\frac{(n-1)}{2} =n(n-1)$. We denote $\overrightarrow{a} = [S_1,S_2, \ldots, S_{\frac{n-1}{2}}]$ such that the stretch $S_k = [(0,2k-1),(2k-1,0), \ldots, (0,2k-1),(2k,0)] = [a_{2n(k-1)+1}, \ldots, a_{2n(k-1)+2n}]$ for all $1 \leq k \leq \frac{n-1}{2}$. 
We denote the $g$-th element of a stretch $S_k$ by $a_{k,g}$ (where $1 \leq g \leq 2n$ and $1 \leq k \leq \frac{n-1}{2}$). In other words  $a_{k,g} = a_{2n(k-1)+g}$ and $S_k =[a_{k,1},a_{k,2},\ldots, a_{k,2n}]$. 

We can remove a sequence of consecutive elements from the left end and/or the right end of a stretch to obtain a (non-empty) \emph{partial stretch}. For any $k$ such that $1 \leq k \leq \frac{n-1}{2}$, we define a partial stretch $S_k^{p,q}$, where  $S_k^{p,q}=[a_{k,p},a_{k,p+1} \ldots ,a_{k,q}]$ where $1\leq p \leq q\leq 2n$. Note that, if $p =1$ and $q = 2n$, then $S_k^{p,q}=S_k$.

The following lemma expresses the total sum over a  partial stretch of the staircase array. Note that all the computations are over $\mathbb{Z}_n$ (that is sums are modulo $n$). We now introduce an additional requirement for $n$ by assuming that $n$ is an odd prime. By imposing this condition, we aim to demonstrate that the total sum of a partial stretch cannot equal $(0,0)$.

\begin{lemma}\label{L: stretch and partial stretches}
Let $n$ be an odd prime and $\overrightarrow{a} = [S_1,S_2, \ldots, S_{\frac{n-1}{2}}]$ be the staircase array. Then for all $S_k \in \overrightarrow{a}$ where  $1 \leq k \leq \frac{n-1}{2}$,  the following hold (where all computations are in $\mathbb{Z}_n$):
\begin{itemize}
    \item [(a)] for a partial stretch $S_k^{p,q}$ where  $1 \leq p \leq q \leq 2n$ we have
        \[\sum_{a_g \in S_k^{p,q}}a_g = 
    \begin{cases}
      ((\lfloor \frac{q}{2}\rfloor-\lfloor \frac{p-1}{2}\rfloor)(2k-1), (\lceil \frac{q}{2} \rceil-\lceil \frac{p-1}{2} \rceil)(2k-1)) & q \not = 2n, \\

      (1-\lfloor \frac{p-1}{2}\rfloor(2k-1), -\lceil \frac{p-1}{2} \rceil(2k-1)) & q  = 2n.
    \end{cases}
    \]

    \item [(b)] $\sum_{a_g \in S_k}a_g = 
    (1,0)$. 

    \item [(c)] $\sum_{a_g \in S_k^{p,q}}a_g \neq (0,0)$.

\end{itemize}
\end{lemma}

\begin{proof}
  \begin{itemize}
    \item [(a)]
   First, assume that $1 \leq p \leq q \leq 2n-1$. Then recalling that $$S_k = [(0,2k-1),(2k-1,0), \ldots, (0,2k-1),(2k,0)] = [a_{2n(k-1)+1}, \ldots, a_{2n(k-1)+2n}]$$ for all $1 \leq k \leq \frac{n-1}{2}$ we have, 

   \begin{align*}
    \sum_{a_g \in S_k^{p,q}}a_g 
    &= \sum_{g=p}^q a_{k,g} \\ 
    &= \sum_{g=1}^q a_{k,g} - \sum_{g=1}^{p-1} a_{k,g} \\ &=   \left(\left\lfloor \dfrac{q}{2} \right\rfloor(2k-1), \left\lceil \dfrac{q}{2} \right\rceil(2k-1) \right) - \left(\left\lfloor \dfrac{p-1}{2} \right\rfloor(2k-1), \left\lceil \dfrac{p-1}{2} \right\rceil(2k-1) \right) \\ 
    &= \left(\left(\left\lfloor \dfrac{q}{2} \right\rfloor- \left\lfloor \dfrac{p-1}{2} \right\rfloor \right)(2k-1), \left(\left\lceil \dfrac{q}{2} \right\rceil-\left\lceil \dfrac{p-1}{2} \right\rceil\right)(2k-1) \right).
   \end{align*}

Now suppose that $q = 2n$. The formula holds if $p=2n$, so assume $p<2n$. Then, 

\begin{align*}
    \sum_{a_g \in S_k^{p,2n}}a_g&=\left(\sum_{a_g \in S_k^{p,2n-1}}a_g \right)+a_{k,2n}\\
    & = \left(\left(n-1- \left\lfloor \dfrac{p-1}{2} \right\rfloor \right)(2k-1), \left(n-\left\lceil \dfrac{p-1}{2} \right\rceil\right)(2k-1) \right)+(2k,0)\\ 
     & = \left(\left(n- \left\lfloor \dfrac{p-1}{2} \right\rfloor \right)(2k-1)+1, \left(n-\left\lceil \dfrac{p-1}{2} \right\rceil\right)(2k-1) \right)\\ 
    &= \left(1-\left\lfloor \frac{p-1}{2}\right\rfloor(2k-1), -\left\lceil \frac{p-1}{2} \right\rceil(2k-1)\right)
   \end{align*}

\item [(b)] Part (b) follows immediately, taking $p=1$ in the previous formula when $q = 2n$.

\item [(c)] Assume $\sum_{a_g \in S_k^{p,q}}a_g = (0,0)$ so that $n$ divides both coordinates of $\sum_{a_g \in S_k^{p,q}}a_g$.
Note that  $n$ is an odd prime and $1 \leq 2k-1 \leq n-2$, so $n$ is coprime with $2k-1$. 

Looking at the second coordinate, we deduce that prime $n$ must divide $\left\lceil \dfrac{q}{2} \right\rceil-\left\lceil \dfrac{p-1}{2} \right\rceil$. Since $1\leq\left\lceil \dfrac{q}{2} \right\rceil\leq n$ and $0\leq \left\lceil \dfrac{p-1}{2} \right\rceil\leq n$, we get that either $p=1$ and $q\in \{2n-1,2n\}$ or $\left\lceil \dfrac{q}{2} \right\rceil=\left\lceil \dfrac{p-1}{2} \right\rceil$. 
The former cannot happen since $\sum_{a_g \in S_k^{1,2n}}a_g = (1,0)$ and $\sum_{a_g \in S_k^{1,2n-1}}a_g = (1-2k,0)$. Thus the latter holds. Since $p\leq q$, this implies that $p=q$ is even. Therefore $\sum_{a_g \in S_k^{p,q}}a_g$ is a single element of $S_k$, none of which is $(0,0)$, a contradiction.

\end{itemize}
\end{proof}

The rest of the section is devoted to proving Theorem~\ref{T: main theorem} by utilizing Lemma~\ref{L: walk is a path} and Lemma~\ref{L: edge orbit general} on the staircase array. This will allow us to derive Proposition~\ref{P: path condition for staircase array}, which provides the necessary and sufficient conditions for a walk to be a path, as well as Proposition~\ref{P: orbit condition for staircase}, which outlines the restrictions on edge orbits.

\begin{proposition} \label{P: path condition for staircase array}
Let $n$ be an odd prime and $\overrightarrow{a}$ be the staircase array. Then $W = W(v_0, \overrightarrow{a})$ is an $n(n-1)$-path.

\end{proposition}

\begin{proof}
Recall that $|\overrightarrow{a}| = n(n-1)$, thus $W(v_0, \overrightarrow{a})$ is an $n(n-1)$-walk. We will show that  $\sum_{g= i+1}^j a_g \neq (0,0)$ for all $0 \leq i< j \leq n(n-1)$. Then by Lemma~\ref{L: walk is a path}, we can conclude that $W(v_0, \overrightarrow{a})$ is a path.

Suppose that $a_{i+1} \in S_k$ and $a_{j} \in S_{k'}$ for some $S_k, S_{k'} \in \overrightarrow{a}$ with $k \leq k'$.
 Thus $a_{i+1}=a_{k,p}= a_{2n(k-1)+p}$ for some $1\leq p\leq 2n$ and $a_{j}= a_{{k'},q}= a_{2n(k'-1)+q}$ for some $1\leq q\leq 2n$.

 \textbf{Case 1:} $k =k'$.

 In this case, there exist  $1 \leq p \leq q \leq 2n$ such that $\sum_{g= i+1}^j a_g = \sum_{a_g \in S_k^{p,q}} a_g \neq (0,0)$  by   Lemma~\ref{L: stretch and partial stretches}(c).

\textbf{Case 2:} $k<k'$.

 This implies $\overrightarrow{a}$ has at least two stretches and, hence $n \geq 5$, since $n$ is an odd prime. Let $k' = k+1+m$ for some $m$ such that $0 \leq m \leq \frac{1}{2}(n-1)-k-1$. This implies that there are $m$ stretches properly between $S_k$ and $S_{k'}$. Therefore there exist  $1 \leq p  \leq 2n$ and  $1 \leq q \leq 2n$ such that $$\sum_{g= i+1}^j a_g =  \sum_{a_g \in S_k^{p,2n}}a_g+m(1,0)+  \sum_{a_g \in S_{k'}^{1,q}}a_g$$ noting that the total sum over a stretch is equal to $(1,0)$ by  Lemma~\ref{L: stretch and partial stretches}(b).

 Now by  Lemma~\ref{L: stretch and partial stretches}(a),  we have:

 \begin{equation} \label{Eq: for i+1 }
  \sum_{a_g \in S_k^{p,2n}}a_g =  \left(1-\left\lfloor \frac{p-1}2 \right\rfloor(2k-1), -\left\lceil \frac{p-1}2 \right\rceil(2k-1)\right)   
 \end{equation}
 
\textbf{Case 2(a):} Suppose that $q <2n$, that is $a_j$ is not the last element of stretch $S_{k'}$.

Then by  Lemma~\ref{L: stretch and partial stretches}(a),  we have:

\begin{equation} \label{Eq: for j }
  \sum_{a_g \in S_{k'}^{1,q}}a_g=  \left( \left\lfloor \frac{q}2 \right\rfloor(2k'-1), \left\lceil \frac{q}2 \right\rceil(2k'-1) \right)   
 \end{equation}

By Equation~\eqref{Eq: for i+1 } and Equation~\eqref{Eq: for j }, we have

\begin{align*}
    \sum_{g=i+1}^{j}a_{g} &= \left(1-\left\lfloor \dfrac{p-1}{2}\right\rfloor(2k-1), -\left\lceil \dfrac{p-1}{2} \right\rceil(2k-1)\right) + m(1,0) + \left(\left\lfloor \dfrac{q}{2}\right\rfloor(2k'-1), \left\lceil \dfrac{q}{2} \right\rceil(2k'-1)\right)\\ 
    &=    \left(m+1-\left\lfloor \dfrac{p-1}{2}\right\rfloor(2k-1)+\left\lfloor \dfrac{q}{2}\right\rfloor(2k'-1), -\left\lceil \dfrac{p-1}{2} \right\rceil(2k-1)+\left\lceil \dfrac{q}{2} \right\rceil(2k'-1)\right)
\end{align*}

If $-\left\lceil \dfrac{p-1}{2} \right\rceil(2k-1)+\left\lceil \dfrac{q}{2} \right\rceil(2k'-1) \ne 0$, then $\sum_{g=i+1}^{j}a_{g} \neq (0,0)$. Therefore, suppose that $-\left\lceil \dfrac{p-1}{2} \right\rceil(2k-1)+\left\lceil \dfrac{q}{2} \right\rceil(2k'-1)  =0$ (recall calculations are in $\mathbb{Z}_n$).
Then noting that $k' = k+1+m$ we have the following:

\[m+1-\left\lfloor \dfrac{p-1}{2}\right\rfloor(2k-1)+\left\lfloor \dfrac{q}{2}\right\rfloor(2k'-1) = 
\begin{cases}
m+1 & \text{if $\lfloor \frac{p-1}{2}
\rfloor = \lceil \frac{p-1}{2} \rceil  \text{ and } \lfloor \frac{q}{2}
\rfloor = \lceil \frac{q}{2} \rceil$},\\
-(m+1) & \text{if $\lfloor \frac{p-1}{2}
\rfloor = \lceil \frac{p-1}{2} \rceil-1 \text{ and } \lfloor \frac{q}{2}
\rfloor = \lceil \frac{q}{2} \rceil-1$},\\
-(2k+m) & \text{if $\lfloor \frac{p-1}{2}
\rfloor = \lceil \frac{p-1}{2} \rceil  \text{ and } \lfloor \frac{q}{2}
\rfloor = \lceil \frac{q}{2} \rceil-1$},\\
2k+m & \text{if $\lfloor \frac{p-1}{2}
\rfloor = \lceil \frac{p-1}{2} \rceil-1  \text{ and } \lfloor \frac{q}{2}
\rfloor = \lceil \frac{q}{2} \rceil$}.
\end{cases}
\]

Note that $m+1 \ne 0 $ because $1 \leq m+1\leq \frac{n-1}{2}-k$, $1 \leq k < \frac{n-1}{2}$ and $5 \leq n$. Therefore,  $-(m+1) \ne 0$ as well. Moreover, $2k+m \ne 0$ because $2k \leq m+2k \leq \frac{n-1}{2}+k-1$ and $1 \leq k < \frac{n-1}{2}$, and hence $-(2k+m) \ne 0 $ as well. 
Therefore, 
\begin{equation*}
    m+1-\left\lfloor \dfrac{p-1}{2}\right\rfloor(2k-1)+\left\lfloor \dfrac{q}{2}\right\rfloor(2k'-1) \ne 0
\end{equation*}
Hence, $\sum_{g=i+1}^{j}a_{g} \neq (0,0)$. 

\textbf{Case 2(b):} Suppose that $q = 2n$, that is $a_j$ is the last element of stretch $S_{k'}$.

Then, $\sum_{a_g \in S_{k'}^{1,q}}a_g=\sum_{a_g \in S_{k'}^{1,2n}}a_g = (1,0)$. Hence,
 \begin{align*}
    \sum_{g=i+1}^{j}a_{g} &= \left(1-\left\lfloor \dfrac{p-1}{2}\right\rfloor(2k-1), -\left\lceil \dfrac{p-1}{2} \right\rceil(2k-1)\right) + m(1,0) + (1,0)\\ &=
    \left(m+2-\left\lfloor \dfrac{p-1}{2}\right\rfloor(2k-1), -\left\lceil \dfrac{p-1}{2} \right\rceil(2k-1)\right)
\end{align*}

If $-\left\lceil \frac{p-1}{2} \right\rceil(2k-1) \ne 0 $, then $\sum_{g=i+1}^{j}a_{g} \neq (0,0)$. Therefore, suppose that $-\left\lceil \frac{p-1}{2} \right\rceil(2k-1) = 0 $.
As noted in the proof of Lemma~\ref{L: stretch and partial stretches}(c), $2k-1$ is coprime to $n$, so  $\left\lceil \frac{p-1}{2} \right\rceil = 0 $ where $1\leq p\leq 2n$. Thus $p=1$ or $p=2n$.
Then,

\[m+2-\left\lfloor \dfrac{p-1}{2}\right\rfloor(2k-1) = 
\begin{cases}
m+2 & \text{if }p=1\\
m+1+2k & \text{if }p=2n

\end{cases}
\]

Note that $m+2 \ne 0 $ because $2 \leq m+2\leq \frac{n-1}{2}-k+1$, $1 \leq k < \frac{n-1}{2}$ and $5 \leq n$. Moreover, $m+1+2k \ne 0 $ because $2k+1 \leq 2k+m+1\leq \frac{n-1}{2}+k$, $1 \leq k < \frac{n-1}{2}$ and $5 \leq n$. Therefore, $\sum_{g=i+1}^{j}a_{g} \neq (0,0)$.

\end{proof}

\begin{proposition}\label{P: orbit condition for staircase}
 Let $n$ be an odd prime and $\overrightarrow{a}$ be the staircase array. Then the path $P = P(v_0, \overrightarrow{a})$ contains at most one edge from each edge orbit under group $G$. 
 
\end{proposition}

\begin{proof}
Let $i,j$ such that $0 \leq i <j < n(n-1)$. We will show that the following holds:
\begin{enumerate} [label=(\alph*)]
    \item  $(\sum_{g=i+1}^j a_g)_2 \neq 0$ or $a_{i+1} \neq a_{j+1}$ and 
    \item $(\sum_{g=i+1}^{j+1} a_g)_2 \neq 0$ or $a_{i+1} \neq -a_{j+1}$.
\end{enumerate}
Then, by Lemma~\ref{L: edge orbit general}, we can conclude that $P$  contains at most one edge from each edge orbit. 

 We first prove (a). If $a_{i+1} \neq a_{j+1}$, then (a) holds. Suppose that $a_{i+1} = a_{j+1}$. Then, due to the nature of a staircase array, both $a_{i+1}$ and $a_{j+1}$ lie in the same stretch. Let $a_{i+1},a_{j+1} \in S_k$, where $k \in \{1,2,\ldots, \frac{n-1}{2}\}$. 
  Recall that
\begin{equation*} 
\begin{aligned} 
    S_k  &= [(0,2k-1),(2k-1,0),\ldots,(2k-1,0),(0,2k-1),(2k,0)] \\ &= [a_{k,1},a_{k,2},\ldots, a_{k,2n}]
\end{aligned} 
\end{equation*}
 
 Then we have $a_{i+1}=a_{k,p}$ and  $a_{j+1}= a_{{k},q}$ where $1 \leq p < q\leq 2n-1$ and $p,q$ have the same parity. In particular $p+2 \leq q$.   By Lemma~\ref{L: stretch and partial stretches}(a) we have,

\begin{equation*}
     \sum_{g=i+1}^j a_g = \sum_{g \in S_k^{p,q-1}} a_g =  \left(\left(\left\lfloor \dfrac{q-1}{2} \right\rfloor- \left\lfloor \dfrac{p-1}{2} \right\rfloor \right)(2k-1), \left(\left\lceil \dfrac{q-1}{2} \right\rceil-\left\lceil \dfrac{p-1}{2} \right\rceil\right)(2k-1) \right).
 \end{equation*}

Then $(\sum_{g=i+1}^j a_g)_2 = \left(\left \lceil \frac{q-1}{2} \right\rceil - \left \lceil \frac{p-1}{2} \right\rceil \right)(2k-1) \neq 0$ since $0 \leq \left \lceil \frac{p-1}{2} \right\rceil < \left \lceil \frac{q-1}{2} \right\rceil \leq n-1$ and $2k-1$ is coprime with $n$. Thus (a) holds.

We now prove (b). If $a_{i+1} \neq -a_{j+1}$, then (b) holds. Suppose that $a_{i+1} = -a_{j+1}$, which is $a_{i+1}+a_{j+1} = (0,0)$. Let $a_{i+1} \in S_k$ and $a_{j+1} \in S_{k'}$, where $k, k' \in \{1,2,\ldots, \frac{n-1}{2}\}$ and $k\leq k'$.  Then 
we have $a_{i+1}=a_{k,p}$ and  $a_{j+1}= a_{{k'},q}$ where $1 \leq p \leq 2n$,  $1 \leq q \leq 2n$ and $p,q$ have the same parity. 

 We now consider two cases according to $a_{i+1},a_{j+1}$ being in the same stretch or two different stretches.
 
\textbf{Case 1:} $k=k'$.

Since $a_{i+1}+a_{j+1} = (0,0)$ and $n$ is an odd prime, both $a_{i+1}$ and $a_{j+1}$ cannot be $(2k-1,0)$ or $(0,2k-1)$ because $1 \leq 2k-1 \leq n-2$ and $n$ is coprime with $2k-1$. Therefore, $a_{j+1} =a_{k,2n} = (2k,0)$ (the last element of $S_k$) and  $a_{i+1} = a_{k,p}=(2k-1,0)$ where $p$ is even. In particular, $2 \leq p \leq 2n-2$.   Now, by Lemma~\ref{L: stretch and partial stretches}(a) we have,

\begin{equation*}
     \sum_{g=i+1}^{j+1} a_g =\sum_{g \in S_k^{p,2n}} a_g=  \left(1-\left\lfloor \frac{p-1}{2} \right\rfloor(2k-1), - \left\lceil \frac{p-1}{2} \right\rceil(2k-1)\right)
 \end{equation*}

Then $(\sum_{g=i+1}^{j+1} a_g)_2 = - \left\lceil \frac{p-1}{2} \right\rceil(2k-1) \neq 0$ since $1 \leq \left\lceil \frac{p-1}{2} \right\rceil \leq n-1$, and $2k-1$ is coprime with $n$. Thus (b) holds if $k=k'$.

\textbf{Case 2:} $k < k'$.

This implies $\overrightarrow{a}$ has at least two stretches and, hence $n \geq 5$, since $n$ is an odd prime. Let $k' = k+1+m$ for some $m$ such that $0 \leq m \leq \frac{1}{2}(n-1)-k-1$. This implies that there are $m$ stretches properly between $S_k$ and $S_{k'}$.
Firstly, observe that $a_{i+1} \in \{(0,2k-1),(2k-1,0),(2k,0)\}$ and $a_{j+1} \in \{(0,2k'-1),(2k'-1,0),(2k',0)\}$. We claim  that  either $a_{i+1} = (2k,0)$ and $a_{j+1} = (2k'-1,0)$, or $a_{i+1} = (2k-1,0)$ and $a_{j+1} = (2k',0)$.

Suppose first that $a_{i+1} = (0,2k-1)$.
If $a_{j+1} \in \{(2k'-1,0),(2k',0)\}$, then $a_{i+1}+a_{j+1} \neq (0,0)$ since  $1 \leq 2k-1 < n-2$. So $a_{j+1}=(0,2k'-1)$, therefore $a_{i+1}+a_{j+1}=(0,2k+2k'-2)$. Note that $2 < 2k+2k'-2 < 2n-4$. Moreover $2k+2k'-2\neq n$ since $n$ is odd. Thus $a_{i+1}+a_{j+1} \neq (0,0)$, a contradiction.

Now suppose that $a_{i+1} = (2k-1,0)$. 
If  $a_{j+1}=(0,2k'-1)$, then $a_{i+1}+a_{j+1} \neq (0,0)$ since  $1 \leq 2k-1 < n-2$. Assume $a_{j+1}=(2k'-1,0)$, therefore $a_{i+1}+a_{j+1}=(2k+2k'-2,0)$ which is not $(0,0)$ by the same argument as above. Thus $a_{i+1}+a_{j+1} \neq (0,0)$ unless $a_{j+1}=(2k',0)$.
 
Lastly, consider when $a_{i+1} = (2k,0)$. 
If  $a_{j+1}=(0,2k'-1)$, then $a_{i+1}+a_{j+1} \neq (0,0)$ since  $1 \leq 2k < n-1$. Assume $a_{j+1}=(2k',0)$, so $a_{i+1}+a_{j+1}=(2k+2k',0)$ which is not $(0,0)$ since  $2 < 2k+2k' < 2n-2$ and $n$ is odd. Thus $a_{i+1}+a_{j+1} \neq (0,0)$ unless $a_{j+1}=(2k'-1,0)$.
This proves the claim.

First, assume that $a_{i+1} = (2k,0) = a_{k,2n}$ and $a_{j+1} = (2k'-1,0) = a_{{k'},q}$ where $q$ is even and $2 \leq q \leq 2n-2$. Then by Lemma~\ref{L: stretch and partial stretches} we have,

\begin{align*}
    \sum_{g=i+1}^{j+1}a_g &= (2k,0)+m(1,0)+\sum_{g \in S_{k'}^{1,q}}a_g\\ 
    &=(2k,0)+m(1,0)+\left(\left\lfloor \dfrac{q}{2}\right\rfloor(2k'-1), \left\lceil \dfrac{q}{2} \right\rceil(2k'-1)\right)\\
    &=\left(2k+m+ \dfrac{q}{2}(2k'-1),  \dfrac{q}{2} (2k'-1)\right).
\end{align*}
Therefore, $(\sum_{g=i+1}^{j+1}a_g)_2 = \dfrac{q}{2} (2k'-1) \neq 0$, because $1 \leq \dfrac{q}{2} \leq n-1$ and  $2k'-1$ is coprime with $n$.

Finally suppose $a_{i+1} = (2k-1,0) = a_{k,p}$ where $p$ is even and $2 \leq p \leq 2n-2 $ and $a_{j+1} = (2k',0) = a_{{k'},2n}$.  Then by Lemma~\ref{L: stretch and partial stretches} we have,

\begin{align*}
    \sum_{g=i+1}^{j+1}a_g &= \sum_{g \in S_{k}^{p,2n}}a_g+m(1,0)+(1,0)\\
    &= \left(1-\left\lfloor \dfrac{p-1}{2} \right\rfloor(2k-1), -\left\lceil \dfrac{p-1}{2} \right\rceil(2k-1)\right) +m(1,0)+(1,0)\\
     &= \left(m+2- \dfrac{p-2}{2}(2k-1), - \dfrac{p}{2} (2k-1)\right).
\end{align*}

Therefore,  $(\sum_{g=i+1}^{j+1}a_g)_2  = -\dfrac{p}{2}(2k-1)\neq 0$, because $1 \leq \dfrac{p}{2}\leq n-1$ and  $2k'-1$ is coprime with $n$. 

Thus (b) holds if $k<k'$.

\end{proof}

Now we have all the necessary steps to prove our main result. Recall that $G = \langle c \rangle$ where $c(a,b) = (a+1,b)$ for all $(a,b) \in \mathbb{Z}_n \times \mathbb{Z}_n$.

\addtocounter{theorem}{-13}

\begin{theorem}

A $G$-transitive $n(n-1)$-path decomposition of $K_n \Box K_n$ exists for all odd primes $n$.

\end{theorem}

\begin{proof}
Let  $n$ be an odd prime, $v_0 = (0,0)$ and $\overrightarrow{a}$ be the staircase array. Then $P = P(v_0, \overrightarrow{a})$ represents an $n(n-1)$-path in $K_n \Box K_n$ by Proposition~\ref{P: path condition for staircase array}. Due to Proposition~\ref{P: orbit condition for staircase} $P$ contains at most one edge from each edge orbit under $G$. Lemma~\ref{L: number of edge orbits} implies that there are $2 \binom{n}{2} = n(n-1)$ edge orbits under $G$. Therefore, $P$ contains exactly one edge from each edge orbit, since $P$ is an $n(n-1)$-path. By Corollary~\ref{C: semiregular} we have that $G$ is semiregular on $E(K_n \Box K_n)$. Therefore, by Theorem~\ref{T: transitive decomposition} $P^G$ is a $G$-transitive $n(n-1)$-path decomposition of $K_n \Box K_n$.

\end{proof}

\section{Conclusion}
We obtained path decompositions of $K_n \Box K_n$, where the paths are arbitrarily long. We also gave some general results on transitive decompositions when the group action is semiregular.
The work here leaves many avenues for further investigation. Our construction involved the staircase array, which requires the primality of $n$. However, it might be possible to look for other constructions which do not require primality.
It would be interesting to establish results similar to ours for all $n$ rather than simply for odd prime $n$.

In this paper we were mainly focusing on the group $G = \langle c \rangle$ where $c(a,b) = (a+1, b)$ for all $(a,b) \in \mathbb{Z}_n \times \mathbb{Z}_m$. In this case, we need $n$ to be odd to maintain the semiregular action on the edges.
We also found $L$-transitive $n(n-1)$-path decompositions of $\mathbb{Z}_n \times \mathbb{Z}_n$ when $n = 2,3,4$, for $L = \langle c' \rangle$ where $c'(a,b) = (a+1, b+1)$ for all $(a,b) \in \mathbb{Z}_n \times \mathbb{Z}_n$. For example, for $n=4$, we take $H = P((0,0), \overrightarrow{a})$ where $\overrightarrow{a} = [(0,1),(0,1),(0,1),(1,0),(1,0),(1,0),(0,2),(2,0),(0,3),(2,0),(0,2),(3,0)]$ and decomposition $H^L$. One potential direction would be to generalize these examples.

Other interesting questions would be considering  $K_n \Box K_m$  for distinct $n$ and $m$, and various possible subgraphs $H$, not just a path. Theorem~\ref{T: transitive decomposition} leaves many possibilities to investigate different $\Gamma$ and $H$.

\bigskip
\noindent\textbf{Acknowledgment.} 
The authors would like to thank Emeritus Professor Cheryl Praeger for fruitful discussions.



\end{document}